\documentclass{amsart}
\usepackage{graphicx}
\usepackage[all]{xy}
\usepackage{amscd}
\usepackage{amssymb}
\usepackage{amsmath}
\usepackage{amsthm}
\usepackage{amsfonts}
\usepackage{graphics}
\usepackage{epsfig,psfrag}
\usepackage[english]{babel}
\usepackage{latexsym}
\usepackage{mathrsfs}

\newtheorem{theorem}{Theorem}[section]
\newtheorem{corollary}[theorem]{Corollary}
\newtheorem{lemma}[theorem]{Lemma}

\newtheorem{proposition}[theorem]{Proposition}
\newtheorem{question}[theorem]{Question}

\theoremstyle{definition}

\newtheorem{definition}[theorem]{Definition}

\theoremstyle{remark}
\newtheorem{remark}[theorem]{Remark}

\newcommand{\RR}{{\mathbb R}}

\newcommand{\ZZ}{{\mathbb Z}}

\begin{document}
\title[Embedding periodic maps of surfaces into those of spheres ]
{Embedding periodic maps of surfaces into those of spheres   with minimal dimensions}


\author{Chao Wang}
\address{School of Mathematical Sciences \& Shanghai Key Laboratory of PMMP, East China Normal University, Shanghai 200241, CHINA}
\email{chao\_{}wang\_{}1987@126.com}

\author{Shicheng Wang}
\address{School of Mathematical Sciences, Peking University, Beijing 100871, CHINA}
\email{wangsc@math.pku.edu.cn}

\author{Zhongzi Wang}
\address{School of Mathematical Sciences, Peking University, Beijing 100871, CHINA}
\email{wangzz22@stu.pku.edu.cn}



\thanks{The first author is supported by National Natural Science Foundation of China (NSFC), grant No. 12131009 and No. 12371067, and Science and Technology Commission of Shanghai Municipality (STCSM), grant No. 18dz2271000.}

\begin{abstract}
It is known that any periodic map of order $n$ on a closed oriented surface of genus $g$  can be equivariantly embedded into $S^m$ for some $m$.  In the orientable and smooth category,  we determine the smallest possible $m$ when  $n\geq 3g$.  We show that for each integer $k>1$ there exist infinitely many periodic maps such that the smallest possible $m$ is equal to $k$.
\end{abstract}

\date{}
\maketitle
\tableofcontents

\section{Introduction}\label{sec:intro}
Suppose  that $M$ is a compact and connected polyhedron  and $G$ is a finite group which acts faithfully on $M$. 
We call such $M$ and $G$ a pair $(M,G)$. 
We  use $\RR^m$,  $S^{m-1}$, $SO(m)$ to denote 
the $m$-dimensional Euclidean space, its unit sphere, its special orthogonal group respectively.
Recall several classical  topics in mathematics.

(1) Embedding $M$ into the $\RR^m$ or $S^m$ is a  topic in topology.
It is known that such embeddings exist 
and to search the minimal $m$ to embed $M$  is an attractive topic.

(2) Embedding finite group $G$ into $SO(m)$ is a topic in algebra. In general
for a given $G$ to find the minimal $m$ is  interesting.

(3) Finite group actions on $M$ is a topic connecting topology and algebra, which is  very active, in particular when $M$ are graphs and surfaces. 

Those topics support a topic about embedding the symmetries of $M$ into the symmetries of $\RR^m$ or $S^m$.
Below we focus on the case of $S^m$, which is more convenient in many cases.
Let  $\text{Aut}(S^m)$ be the group of automorphisms of $S^m$, where automorphisms  can be either homeomorphisms,
or diffeomorphisms, or isometries. 

\begin{definition} For a given pair $(M, G)$,
call an embedding
$e: M\to S^m$  $G$-equivariant, 
if  there is a representation $ \rho: G\to Aut(S^n)$ such that 
$$e\circ g= \rho(g) \circ e$$
for any $g\in G$.
\end{definition}


The existence of $G$-equivariant embedding for any pair $(M, G)$
follows from the classical work of Mostow \cite{Mos} and  also Palais \cite{Pal}.

\begin{question}\label{main} For a given pair $(M, G)$, how to find  the minimal $m$, so that there is a $G$-equivariant embedding $M\to S^m$.
\end{question}

 For a given pair $(M, G)$, some upper  and lower bounds of $m$ in Question \ref{main} are given in \cite{Was}, \cite{Wz2}. 
But those  bounds are far away from sharp in general.

There are systematic
studies on  $G$-equivariant embeddings for graphs and surfaces into $S^3$, see \cite{Cos}, \cite{FNPT}, \cite{WWZZ}, \cite{NWW} and the references therein. Those studies rely on the geometry and topology of 3-manifolds developed in the last several decades, as well as  our 3-dimensional intuition.  Once we know that there is a $G$-equivariant embedding  $M\to S^3$, then usually $3$ is  the minimum
$m$ for the pair  $(M, G)$ in Question \ref{main}, since usually graphs and surfaces themselves can not be embedded into $S^2$.
  
If there is no  $G$-equivariant embedding $M\to S^3$ for a  pair $(M,G)$,
then to determine the minimal $m$ in Question \ref{main} is  a challenge.
 Results   for $m>3$ in Question \ref{main} appear recently for some interesting pairs, see \cite{Zim}, \cite {Wz1}.  

In this paper, we will focus on a primary case: the finite cyclic group $\ZZ_n$ actions on $\Sigma_g$, the closed orientable surfaces of genus $g$.
We will work in the {\it smooth and orientable} category (unless otherwise stated), i.e., in Definition 1.1., the group actions on $\Sigma_g$ and on $S^n$ are smooth and orientation-preserving,
and  embeddings $\Sigma_g\to S^m$ are smooth.
Now  each $\ZZ_n$ action on $\Sigma_g$ is determined by an orientation-preserving  smooth periodic map $f$, will be denoted as $(\Sigma_g, f)$. For such $(\Sigma_g, f)$, we use $D_g(f)$ to denote $m$ in Question \ref{main},
and call $D_g(f)$ the minimal embedding dimension for $(\Sigma_g, f)$. 



\begin{question}\label{que:Dgf}
For a periodic map $f$ on $\Sigma_g$, what is the exact value of $D_g(f)$?
\end{question}

Let $f$ be a periodic map  on $\Sigma_g$ of order $n$.
Clearly $D_0(f)=2$ and $D_g(f)\geq 3$ when $g\geq 1$.  Call $f$ is free if each $f$-orbit has length  $n$.
An   old result of Nielsen claims that  $D_g(f)=3$ when $f$ is free \cite{Ni}, and in this case the $n\le g-1$ unless $g\le 1$.
Furthermore, the periodic maps with $D_g(f)=3$ are classified recently in \cite{NWW}.  

In this paper, we will determine $D_g(f)$ for
the case of $f$ with order $n \geq 3g$.

\begin{theorem}\label{main1} Let  $f$ be a periodic map  on $\Sigma_g$ with order $n\geq 3g$, $g>0$. 

Then $D_g(f)=6$ with the following exceptions: $D_2(f)= 3$ or $4$, which are determined, if $f$  is of  order $6$;
$D_1(f)=5$ if $f$  is of  order $3$ and is non-free;   $D_1(f)=4$    if $f$  of  order $4$ and is non-free; $D_1(f)=3$ if $f$ is free.
\end{theorem}

\begin{remark} 
Since $D_1(f)=3$ when a periodic map $f$ is of order 2 \cite{NWW}, $D_1(f)$ for all periodic maps $f$ on the torus are determined. 
\end{remark}

Also, we have the following result.

\begin{theorem}\label{main2}
For any integer $m>1$, there exist infinitely many $(\Sigma_g, f)$, where $f$ is a periodic map on $\Sigma_g$, such that $D_g(f)=m$.
\end{theorem}

Furthermore, if we still consider smooth and orientation-preserving actions but do not require the embedding $e$ to be smooth, then we can similarly define another minimum embedding dimension $\hat{D}_g(f)$. Clearly $\hat{D}_g(f)\le {D}_g(f)$. As a comparison, we present 
the following


 

\begin{proposition}\label{thm:s}
$\hat{D}_g(f)=D_g(f)$ for all maps $f$ in the proof of Theorem \ref{main2}.

 $\hat{D}_g(f)<{D}_g(f)$ for  all maps $f$  with  ${D}_g(f)=6$  in Theorem \ref{main1} except one.
\end{proposition}

{\bf About the Proofs.} For each $g>0$, the periodic maps  $f$ having order $n\geq 3g$ are classified in terms of their orders and information of their orbifolds $\Sigma_g/f$ in \cite{Br} for $g\le 3$ and in \cite{Hi} for $g\ge 3$. In Section~\ref{sec:examp},
we will construct an explicit $f$-equivariant embedding  $\Sigma_g\to S^m$ for each involved periodic map $f$, which gives an upper bound of 
$D_g(f)$. In Section~\ref{sec:proof},  by analyzing the tangent spaces at fixed point sets and by applying Smith theory, we get a lower bound of $D_g(f)$, which coincides with the upper bound,  hence finish the proofs of Theorem \ref{main1} and Theorem \ref{main2}. 

The constructions in Section~\ref{sec:examp} have three steps. 
Step 1: Construct embeddings of punctured surfaces into $S^3$ smoothly and equivariantly. For maps involved in Theorem \ref{main1},   the required symmetries of punctured surfaces are directly embedded into the symmetries of $S^3$, which are more intuitive. For maps involved in Theorem \ref{main2},
the constructions start from the embeddings of  compact surfaces into spherical 3-orbifolds, which are more sophisticated.
Step 2: Cap off the boundaries by smooth disks in $S^m\supset S^3$ equivariantly via the geometric join trick. 
Step 3:  Smooth the embeddings. 
The last step is not difficult since we only have one dimensional singularities and on which the action is free.

\section{Constructions of equivariant embeddings}\label{sec:examp}
In this section we will give several constructions of equivariant embeddings. The contents in Section~\ref{subsec:Borsur} and Section~\ref{subsec:Closur} will provide upper bounds of $D_g(f)$ for periodic maps in Theorems~\ref{main1}, and some ideas will also be used in Section~\ref{subsec:orbifold}, which will provide upper bounds of $D_g(f)$ for periodic maps in Theorems~\ref{main2}.

\subsection{Equivariant bordered surfaces in $S^3$}\label{subsec:Borsur}
Write $\mathbb{R}^4=\mathbb{C}^2$,
\[S^3=\{(z,z')\in\mathbb{C}^2\mid |z|^2+|z'|^2=1\}.\]
The Clifford Torus is given by
\[T=\{(z,z')\in S^3\mid |z|=|z'|\}.\]
It splits $S^3$ into two handlebodies
\[H=\{(z,z')\in S^3\mid |z|\geq|z'|\}, H'=\{(z,z')\in S^3\mid |z|\leq|z'|\}.\]
Let $p$ and $q$ be two positive integers. Define disks
\[D_j=\{(z,z')\in H\mid p\arg z=2j\pi\}, D_k'=\{(z,z')\in H'\mid q\arg z'=2k\pi\},\]
where $j=0,1,\ldots,p-1$ and $k=0,1,\ldots,q-1$. Define an annulus
\[A'=\{(z,z')\in H'\mid \arg z=q\arg z'\}.\]
Define $\tau_n$ to be the $2\pi/n$-rotation on the unit circle $S^1\subset \mathbb{C}$.
Define an isometry $\tau_{p,q}:\mathbb{C}^2\rightarrow\mathbb{C}^2$ by
\[(z,z')\mapsto(e^{i\frac{2\pi}{p}}z,e^{i\frac{2\pi}{q}}z'),\]
then $\tau_{p,q}$ keeps $S^3$, $T$, $H$, $H'$, $\{D_j\}$ and $\{D_k'\}$ invariant. The isometry $\tau_{p,pq}$ keeps $\{D_j\}$ and $A'$ invariant.

\begin{figure}[h]
\includegraphics{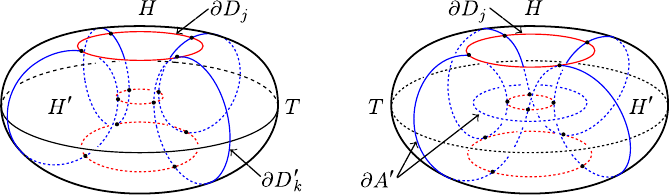}
\caption{The boundaries of $D_j$, $D_k'$ and $A'$}\label{fig:DiskAnnulus}
\end{figure}

Let $\mathcal{P}:S^3\setminus\{(-1,0)\}\rightarrow\mathbb{R}^3$ be the stereographic projection given by
\[(u+vi,x+yi)\mapsto(\frac{v}{1+u},\frac{x}{1+u},\frac{y}{1+u}),\]
then $\mathcal{P}(T)$, $\mathcal{P}(H)$ and $\mathcal{P}(H')$ will be like in Figure~\ref{fig:DiskAnnulus}. The left one shows $\mathcal{P}(\partial D_j)$ and $\mathcal{P}(\partial D_k')$ while the right one shows $\mathcal{P}(\partial D_j)$ and $\mathcal{P}(\partial A')$, where $p=3,q=4$. Note that under the spherical metric all $\partial D_j$ and $\partial D_k'$ have the same length and the component of $\partial A'$ in the interior of $H'$ is the geometric core of $H'$.

\begin{figure}[h]
\includegraphics{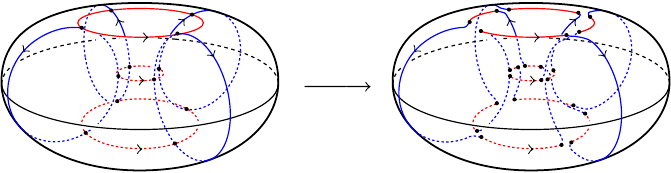}
\caption{Surgeries on the boundaries of $D_j$ and $D_k'$}\label{fig:CurveSur}
\end{figure}

The boundaries $\{\partial D_j\}$ and $\{\partial D_k'\}$ intersect at $pq$ points on $T$. We can apply surgeries in $T$ equivariantly at each intersection, then the set $\{\partial D_j\}\cup\{\partial D_k'\}$ will become a torus link as shown in Figure~\ref{fig:CurveSur}.

There are two ways to do the equivariant surgeries and the results are mirror images of each other. We will fix a way as follows: The boundaries $\partial D_j$ and $\partial D_k'$ can be parameterized respectively by
\[(\frac{\sqrt{2}}{2}e^{i\frac{2j\pi}{p}},\frac{\sqrt{2}}{2}e^{i\theta}), (\frac{\sqrt{2}}{2}e^{i\theta},\frac{\sqrt{2}}{2}e^{i\frac{2k\pi}{q}}), \theta\in [0,2\pi],\]
then the $\theta$-increasing direction gives directions of $\partial D_j$ and $\partial D_k'$, and the surgery should keep the directions. Surgeries in Figure~\ref{fig:CurveSur} are applied in this way.

\begin{figure}[h]
\includegraphics{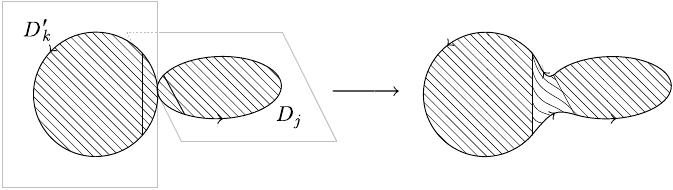}
\caption{Surgeries on the disks $D_j$ and $D_k'$}\label{fig:DiskSur}
\end{figure}

According to the equivariant surgeries on the set $\{\partial D_j\}\cup\{\partial D_k'\}$, we can also apply surgeries on $\{D_j\}\cup\{D_k'\}$ equivariantly. At each intersection, we can choose a small neighborhood and replace the disks in it by a band as shown in Figure~\ref{fig:DiskSur}.

Via the equivariant surgeries $\{D_j\}\cup\{D_k'\}$ will become a connected bordered surface in $S^3$, which is a Seifert surface of a torus link. By the construction, the map $\tau_{p,q}$ keeps the surface invariant. Denote this surface by $S_{p,q}$.

Similarly, the boundaries $\{\partial D_j\}$ and $\partial A'$ intersect at $pq$ points on $T$, and we can apply equivariant surgeries as above. Now the two ways to do the equivariant surgeries may give essentially different results.

Denote the boundary of $A'$ in $T$ by $\partial_1$ and the other boundary by $\partial_o $. Then the boundaries $\partial_o $ and $\partial_1$ can be parameterized respectively by
\[(0,e^{-i\theta}),(\frac{\sqrt{2}}{2}e^{iq\theta},\frac{\sqrt{2}}{2}e^{i\theta}), \theta\in [0,2\pi],\]
and the $\theta$-increasing direction gives directions of $\partial_o $ and $\partial_1$.

By applying equivariant surgeries which keep the directions of $\partial D_j$ and $\partial_1$, we can obtain a connected bordered surface in $S^3$. Denote it by $S^+_{p,q}$. If we reverse the directions of $\{\partial D_j\}$ and apply the equivariant surgeries, then we will get another connected bordered surface in $S^3$. Denote this surface by $S^-_{p,q}$. 

By the construction, the map $\tau_{p,pq}$ keeps $S^+_{p,q}$ and $S^-_{p,q}$ invariant. The boundary of $S^+_{p,q}$ (resp. $S^-_{p,q}$) consists of two parts: $\partial_o $ and a torus link in $T$. Denote the torus link by $\partial^+_{p,q}$ (resp. $\partial^-_{p,q}$). We also use $\partial_{p,q}$ to denote the boundary of $S_{p,q}$.

\begin{figure}[h]
\includegraphics[width=300pt, height=150pt]{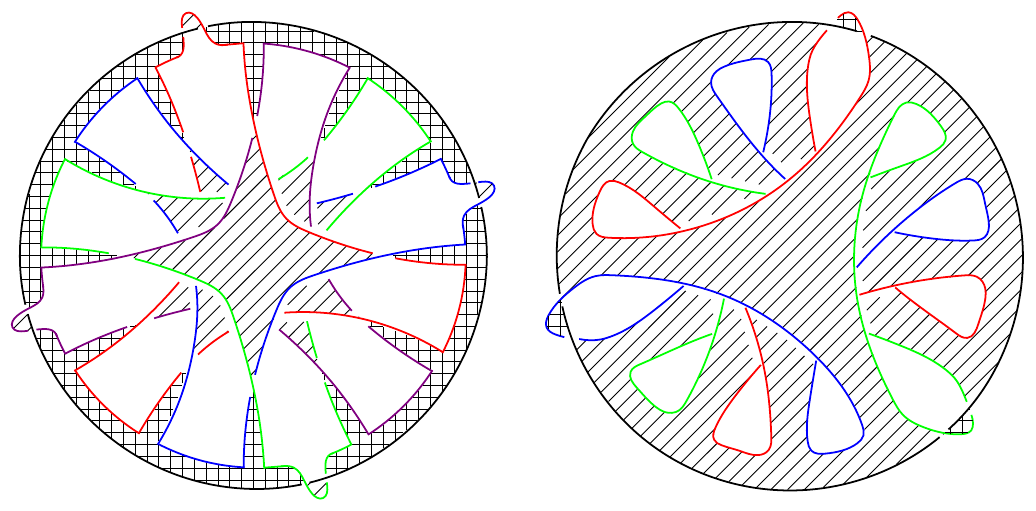}
\caption{$S^+_{3,4}$ with $\#\partial^+_{3,4}=4$, and $S^-_{4,3}$ with $\#\partial S^-_{4,3}=3$}\label{fig:SPMpq}
\end{figure}

In Figure~\ref{fig:SPMpq}, the left one is a sketch of $S^+_{3,4}$, where $\partial S^+_{3,4}$ has 4 components; while the right one is a sketch of $S^-_{4,3}$, where $\partial S^-_{4,3}$ has 3 components.


{\bf Notations.} We will use $(a\,;b)$ and $[a\,;b]$ to denote the greatest common divisor and the least common multiple of integers $a$ and $b$, respectively.

\begin{lemma}\label{lem:toptype}
(1) The surface $S_{p,q}$ is an orientable surface of genus
\[g=\frac{pq-p-q-(p\,;q)}{2}+1,\]
and $\partial_{p,q}$ consists of $d$ torus knots of type $(p/d,q/d)$, where $d=(p\,;q)$.

(2) The surface $S^+_{p,q}$ is an orientable surface of genus
\[g=\frac{pq-p-(p+1\,;q)+1}{2},\]
and $\partial^+_{p,q}$ consists of $d$ torus knots of type $((p+1)/d,q/d)$, where $d=(p+1\,;q)$.

(3) The surface $S^-_{p,q}$ is an orientable surface of genus
\[g=\frac{pq-p-(p-1\,;q)+1}{2},\]
and $\partial^-_{p,q}$ consists of $d$ torus knots of type $((p-1)/d,q/d)$, where $d=(p-1\,;q)$.
\end{lemma}

\begin{proof}
In above constructions, we have defined orientations on $\partial D_j$, $\partial D_k'$ and $\partial_1$. Since the surgeries keep the orientations, the surfaces are all orientable.

In each case the surface contains a graph as its deformation retract, where the surgery band corresponds to an edge. Hence, we can get the Euler characteristic of the surface. Then, we only need to proof the part about the boundaries.

We regard $T$ as $\partial D_0'\times\partial D_0$, which has a natural Euclidean structure such that $\partial D_0'$ and $\partial D_0$ have length $2\pi$ and are orthogonal to each other.

In (1), since $\partial D_j$ has slope $\infty$ and $\partial D_k'$ has slope $0$, the components of $\partial_{p,q}$ will have slope $p/q$. Hence, each component of $\partial_{p,q}$ is a $(p/d,q/d)$-torus knot, and $\partial_{p,q}$ has $d$ components. Since the Euler characteristic
\[\chi(S_{p,q})=p+q-pq,\]
we can get the required formula of $g$.

In (2), since $\partial D_j$ has slope $\infty$ and $\partial_1$ has slope $1/q$, the components of $\partial^+_{p,q}$ will have slope $(1+p)/q$. Hence, each component is a $((p+1)/d,q/d)$-torus knot, and $\partial^+_{p,q}$ has $d$ components. Then $\partial S^+_{p,q}$ has $d+1$ components, and combined with
\[\chi(S^+_{p,q})=p-pq,\]
we can get the required formula of $g$.

The proof of (3) is similar to (2), except that the components of $\partial^-_{p,q}$ will have slope $(1-p)/q$ and each component is a $((p-1)/d,q/d)$-torus knot.
\end{proof}

In each case, the boundaries of the surface are equipped with orientations from the construction. Since $\tau_{p,q}$ and $\tau_{p,pq}$ act as translations on $T$, $T/\tau_{p,q}$ and $T/\tau_{p,pq}$ are tori, with the induced Euclidean structures, and the images of the boundaries in them are isotopic to geodesics. Hence we can assume that $\partial_{p,q}$ (resp. $\partial^+_{p,q}$, $\partial^-_{p,q}$) has a standard parametrization by $[0,2\pi]$, as in Lemma~\ref{lem:trans} below. Then, on each component of the boundaries of the surface, some power of $\tau_{p,q}$ or $\tau_{p,pq}$ will induce a translation along the circle by some $\theta$, which we will call a $\theta$-rotation.

\begin{lemma}\label{lem:trans}
(1) The map $\tau_{p,q}$ acts trivially on the set of components of $\partial_{p,q}$, and it induces a $2\pi/[p\,;q]$-rotation on each component of $\partial_{p,q}$.

(2) The map $\tau_{p,pq}$ acts transitively on the set of components of $\partial^+_{p,q}$ and induces a $-2\pi/(pq)$-rotation on $\partial_o $. Let $d=(p+1\,;q)$ and $s$ be an integer satisfying
\[\frac{1+p}{d}s\equiv 1\,(\mathrm{mod}\,\frac{pq}{d}).\]
Then the map $\tau^d_{p,pq}$ induces a $2ds\pi/(pq)$-rotation on each component of $\partial^+_{p,q}$.

(3) The map $\tau_{p,pq}$ acts transitively on the set of components of $\partial^-_{p,q}$ and induces a $-2\pi/(pq)$-rotation on $\partial_o $. Let $d=(p-1\,;q)$ and $s$ be an integer satisfying
\[\frac{1-p}{d}s\equiv 1\,(\mathrm{mod}\,\frac{pq}{d}).\]
Then the map $\tau^d_{p,pq}$ induces a $2ds\pi/(pq)$-rotation on each component of $\partial^-_{p,q}$.
\end{lemma}

\begin{proof}
In (1), by Lemma~\ref{lem:toptype}(1), the $d$ components of $\partial_{p,q}$ can be given by
\[(\frac{\sqrt{2}}{2}e^{i\frac{q\theta}{d}},\frac{\sqrt{2}}{2}e^{i\frac{p\theta}{d}+i\frac{2l\pi}{q}}), \theta\in [0,2\pi],\]
where $l=0,1,\ldots,d-1$. They can also be parameterized by
\[(\frac{\sqrt{2}}{2}e^{i\frac{q\theta}{d}+i\frac{2bl\pi}{d}},\frac{\sqrt{2}}{2}e^{i\frac{p\theta}{d}+i\frac{2al\pi}{d}}), \theta\in [0,2\pi],\]
where $a$ and $b$ are two integers satisfying $aq-bp=d$. Then, the action of $\tau_{p,q}$ on the components is given by
\[(\frac{\sqrt{2}}{2}e^{i\frac{q\theta}{d}+i\frac{2bl\pi}{d}},\frac{\sqrt{2}}{2}e^{i\frac{p\theta}{d}+i\frac{2al\pi}{d}})\mapsto (\frac{\sqrt{2}}{2}e^{i\frac{q}{d}(\theta+\frac{2d\pi}{pq})+i\frac{2bl\pi}{d}}, \frac{\sqrt{2}}{2}e^{i\frac{p}{d}(\theta+\frac{2d\pi}{pq})+i\frac{2al\pi}{d}}).\]
Hence, on each component $\theta$ increases $2\pi/[p\,;q]$ and we have the result.

In (2), by Lemma~\ref{lem:toptype}(2), the $d$ components of $\partial^+_{p,q}$ can be given by
\[(\frac{\sqrt{2}}{2}e^{i\frac{q\theta}{d}},\frac{\sqrt{2}}{2}e^{i\frac{(1+p)\theta}{d}+i\frac{2l\pi}{q}}), \theta\in [0,2\pi],\]
where $l=0,1,\ldots,d-1$. They can also be parameterized by
\[(\frac{\sqrt{2}}{2}e^{i\frac{q\theta}{d}+i\frac{2bl\pi}{d}},\frac{\sqrt{2}}{2}e^{i\frac{(1+p)\theta}{d}+i\frac{2al\pi}{d}}), \theta\in [0,2\pi],\]
where $a$ and $b$ are two integers satisfying $aq-b(1+p)=d$. The action of $\tau_{p,pq}$ on the components maps the above formula to
\[(\frac{\sqrt{2}}{2}e^{i\frac{q}{d}(\theta+\frac{2(d+bp)\pi}{pq})+i\frac{2b(l-1)\pi}{d}}, \frac{\sqrt{2}}{2}e^{i\frac{1+p}{d}(\theta+\frac{2(d+bp)\pi}{pq})+i\frac{2a(l-1)\pi}{d}}).\]
Hence, $\tau_{p,pq}$ acts transitively on the set of components of $\partial^+_{p,q}$. We can also obtain this directly from the construction, because each component contains arcs from $\partial_1$. According to the orientation of $\partial_o $, the map $\tau_{p,pq}$ induces a $-2\pi/(pq)$-rotation on it. Now we only need to show that $s\equiv d+bp\,(\mathrm{mod}\,pq/d)$. This is because
\[\frac{1+p}{d}\cdot(d+bp)=\frac{(1+p)(aq-b)}{d}=1+a\cdot\frac{pq}{d}\equiv 1\,(\mathrm{mod}\,\frac{pq}{d}).\]
Hence, the map $\tau_{p,pq}^d$ induces a $2ds\pi/(pq)$-rotation on each component of $\partial^+_{p,q}$.

The proof of (3) is similar to (2). The components of $\partial^-_{p,q}$ can be given by
\[(\frac{\sqrt{2}}{2}e^{i\frac{q\theta}{d}+i\frac{2bl\pi}{d}},\frac{\sqrt{2}}{2}e^{i\frac{(1-p)\theta}{d}+i\frac{2al\pi}{d}}), \theta\in [0,2\pi],\]
where $a$ and $b$ are two integers satisfying $aq-b(1-p)=d$, and the action of $\tau_{p,pq}$ on the components maps the above formula to
\[(\frac{\sqrt{2}}{2}e^{i\frac{q}{d}(\theta+\frac{2(d-bp)\pi}{pq})+i\frac{2b(l+1)\pi}{d}}, \frac{\sqrt{2}}{2}e^{i\frac{1-p}{d}(\theta+\frac{2(d-bp)\pi}{pq})+i\frac{2a(l+1)\pi}{d}}).\]
One can also verify that $s\equiv d-bp\,(\mathrm{mod}\,pq/d)$.
\end{proof}

Let $\hat{S}_{p,q}$ be the closed surface obtained by capping each component of $\partial S_{p,q}$ by a disk. Since rotations on the boundary of a disk can extend into its interior, the restriction of $\tau_{p,q}$ on $S_{p,q}$ extends to a periodic map on $\hat{S}_{p,q}$. Denote the extension by $f_{p,q}$. Similarly, we can have $\hat{S}^+_{p,q}$, $\hat{S}^-_{p,q}$, $f^+_{p,q}$ and $f^-_{p,q}$, where $f^+_{p,q}$ and $f^-_{p,q}$ are the periodic extensions of the restrictions of $\tau_{p,pq}$ on $S^+_{p,q}$ and $S^-_{p,q}$, respectively.

Let $f$ be a periodic map on $\Sigma_g$ of order $n$. Let $\Sigma_g/f$ denote the orientable closed 2-orbifold determined by $f$. We say that $\Sigma_g/f$ has {\it type $(\bar{g}:n_1,\ldots,n_l)$} if it has underlying space $\Sigma_{\bar{g}}$ and except $l$ points of indices $n_1,\ldots,n_l$ all its points are regular. Here a point in $\Sigma_g/f$ has {\it index} $d$ if it has $n/d$ preimages in $\Sigma_g$, and it is {\it regular} if it has index $1$, otherwise it is {\it singular}. The type of $\Sigma_g/f$ is unique up to reordering $n_1,\ldots,n_l$ and removing all $n_j=1$.

\begin{lemma}\label{lem:types}
(1) Let $d=(p\,;q)$ and $n_l=[p\,;q]$ for $l=1,\ldots,d$, then the orbifold $\hat{S}_{p,q}/f_{p,q}$ has type $(0:n_1,\ldots,n_d,p/d,q/d)$.

(2) The orbifold $\hat{S}^+_{p,q}/f^+_{p,q}$ has type $(0:pq,pq/d,q)$, where $d=(p+1\,;q)$.

(3) The orbifold $\hat{S}^-_{p,q}/f^-_{p,q}$ has type $(0:pq,pq/d,q)$, where $d=(p-1\,;q)$.
\end{lemma}

\begin{proof}
In (1), by Lemma~\ref{lem:trans}(1), those $d$ disks bounded by the components of $\partial_{p,q}$ give $d$ points of index $[p\,;q]$ in $\hat{S}_{p,q}/f_{p,q}$. Since $\tau_{p,q}$ acts transitively on $\{D_j\}$, these disks give a point of index $q/d$ in $\hat{S}_{p,q}/f_{p,q}$. Similarly, the disks $\{D_k'\}$ give a point of index $p/d$ in $\hat{S}_{p,q}/f_{p,q}$. Hence, $\hat{S}_{p,q}/f_{p,q}$ has the type $(\bar{g}:n_1,\ldots,n_d,p/d,q/d)$ for some $\bar{g}$. By Lemma~\ref{lem:toptype}(1) and the Riemann--Hurwitz formula, we have
\[-pq+p+q+(p\,;q)=[p\,;q](2-2\bar{g}-d\cdot\frac{[p\,;q]-1}{[p\,;q]}-\frac{p-d}{p}-\frac{q-d}{q}).\]
Hence, we have $\bar{g}=0$.

In (2), by Lemma~\ref{lem:trans}(2), the disk bounded by $\partial_o $ gives a point of index $pq$ and those $d$ disks bounded by the components of $\partial^+_{p,q}$ give a point of index $pq/d$. Since $\tau_{p,pq}$ acts transitively on $\{D_j\}$, these disks give a point of index $q$. Hence, $\hat{S}^+_{p,q}/f^+_{p,q}$ has the type $(\bar{g}:pq,pq/d,q)$ for some $\bar{g}$. As in (1), we have
\[-pq+p+(p+1\,;q)+1=pq(2-2\bar{g}-\frac{pq-1}{pq}-\frac{pq-d}{pq}-\frac{q-1}{q}).\]
Hence, we have $\bar{g}=0$.

The proof of (3) is similar to (2).
\end{proof}

A consequence of Lemma \ref{lem:types} (1) is the following fact.

\begin{corollary}\label{later} Let $d=(p\,;q)$. Then the orbifold ${S}_{p,q}/\tau_{p,q}$ is a $S^2$ with $d$ disks removed  and has 
at most two singular
points of indices   $p/d$ and $q/d$.
\end{corollary}

\subsection{Equivariant closed surfaces in $S^m$}\label{subsec:Closur}

In this section we will construct embeddings for Theorem \ref{main1}, and some ideas and results will be used to construct 
embeddings for  Theorem \ref{main2}.

We consider $S^p$ and $S^q$ as  the unit spheres of $\RR^{p+1}$ and $\RR^{q+1}$ respectively, and consider  $\RR^{p+q+2}$ as the product 
$\RR^{p+1}\times \RR^{q+1}$. Here we identify $\mathbb{R}^{p+1}$ with $\mathbb{R}^{p+1}\times \{0\}$ canonically and identify $\mathbb{R}^{q+1}$ with $\{0\}\times \mathbb{R}^{q+1}$ canonically. Then there is a geometrical way to join  a subset $X$ of  $S^p$ and a subset $Y$ of  $S^q$ 
to provide $X*Y\subset S^{p+q+1}\subset \RR^{p+q+2}$ as below: 

For each point $x\in S^p\subset \RR^{p+1}$ and each point $y\in S^q\subset \RR^{q+1}$, the line segment  $[x,y]$ to join $x$ and $y$ 
is given by the following parametrized curve in variable $t$. 
\[ [x,y]_t=x {\rm cos} \frac{\pi t}2+ y {\rm sin} \frac{\pi t}2, \,\, t\in [0,1].\]
Now we define $$X*Y=\bigcup_{x\in X, y\in Y}{[x,y]}.$$
For a point $x$, $x*Y$ is a cone from $x$ over $Y$. 
 One can verify that 

\begin{lemma}\label{join}

(1) $[x, y]$ belongs to $S^{p+q+1}$, indeed is the unique shortest geodesic in $S^{p+q+1}$
joining $x$ and $y$.

(2) The interiors of any two segments are disjoint and each point of $S^{p+q+1}$ lies in some segment $[x,y]$.
Hence $S^p* S^q= S^{p+q+1}$. 

(3) For any $\tau_1 \in O(p+1)$ acting on $S^p$ and $\tau_2 \in O(q+1)$ acting on $S^q$, $\tau_1 \oplus \tau_2\in O(p+q+2)$
acting on $S^{p+q+1}$ by  

$$\tau_1 \oplus \tau_2 (x cos \frac{\pi t}2+ y sin \frac{\pi t}2)=\tau_1  (x) cos \frac{\pi t}2+ \tau_2 (y) sin \frac{\pi t}2.$$

In addition: If $\tau_1(X)=X$  and  $\tau_2(Y)=Y$, then $\tau_1\oplus\tau_2 (X*Y)=X*Y$, where $X\subset S^p, Y\subset S^q$.

(4) Suppose that $x\in S^p$  is a point, and  $C$ is a great circle in $S^q$, then $x*C$ is a geodesic disk bounded by $C$.


\end{lemma}

\begin{proposition}\label{proof}
Let  $\hat S$ be either
$\hat{S}_{p,q}$, or $\hat{S}^{\pm}_{p,q}$; 
and $(f, \hat S)$ be either
$(f_{p,q}, \hat{S}_{p,q})$, or $(f^{\pm}_{p,q}, \hat{S}^{\pm}_{p,q})$.
Then $(f|S, S)$ is either
$(\tau_{p,q}, {S}_{p,q})$, or $(\tau_{p,pq}, {S}^\pm_{p,q})$, and $(f|S, S)\subset (SO(4), S^3)$, 

(a) If either $\partial_{p,q}$ (resp. $\partial^{\pm}_{p,q}$) is connected or each component of $\partial_{p,q}$   (resp. $\partial^{\pm}_{p,q}$) is unknotted in $S^3$, 
then there is an $ f$-equivariant embedding  $e:\hat S\rightarrow S^6$.

(b) In the case of each component of $\partial_{p,q}$ being unknotted, the $S^6$ can be replaced by $S^5$, and furthermore by $S^4$ if $(p,q)=2$.
\end{proposition}
 
\begin{proof} (a) We just prove it for $S^-_{p,q}$, the proof for $S^+_{p,q}$ is the same, and some differences of the proof for $S_{p,q}$ will be illustrated in the proof of (b). 

Now we have the embedding 
$(\tau_{p,pq}, {S}^-_{p,q})\subset (SO(4), S^3).$
Let $N$ be the north pole of $S^2$ and $D^-$ be the down-half sphere of $S^2$ with  the equator $S_E^1$ as the boundary.
Schematic pictures Figures 5 and 6 provide intuition of the proof below.

(i) Suppose $\partial^{-}_{p,q}$ is connected: Then $d=(p-1; q)=1$ and $\partial^{-}_{p,q}$ is a torus knot of type $(p-1, q)$ by Lemma \ref{lem:toptype} (3),
and $\tau_{p,pq}$ induces a $-2\pi/(pq)$-rotation on $\partial_o $, and 
a $2s\pi/(pq)$-rotation on $\partial^ -_{p,q}$ for some $s$ by Lemma \ref{lem:trans} (3).

\begin{figure}[h]
\includegraphics[width=330pt, height=135pt]{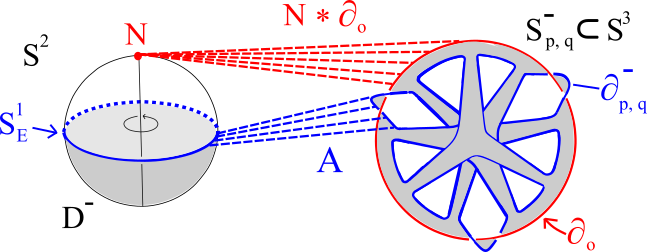}
\caption{Sketch for $f^-_{p,q}$-equivariant embedding $\hat{S}^-_{p,q}\to S^2*S^3=S^6$ when $\partial^-_{p,q}$ is connected}\label{fig:CurveSur1}
\end{figure}

Now consider the join $S^2*S^3=S^6$:
Let $\iota$ be the $2s\pi/(pq)$-rotation of $S^2$ around $N$.
Since $\partial _0$ is a great circle with $\tau_{p,pq}(\partial_o )=\partial_o $, and $\iota(N)=N$, by Lemma \ref{join} (3) and (4),  the join $N*\partial _0$ is a $\iota\oplus \tau_{p,pq}$ invariant  geodesic disc,  which caps on $\partial_o \subset S^-_{p,q}$, where part of $N*\partial _0$ is indicated  by the dashed red lines in Figure 5.
Let $S^1_E(\theta), \partial^-_{p,q}(\theta)$ be smooth equivariant parametrizations (i. e. $\iota\circ S^1_E(\theta) =S^1_E(\theta+2\pi s/pq)$, $\tau_{p.pq}\circ\partial_o (\theta)=\partial_o (\theta+2\pi s/pq)$).
Let $$ A=\bigcup_{\theta\in [0,2\pi]} [S^1_E(\theta), \partial^-_{p,q}(\theta)].$$
Then $A$ is a  smooth annulus connecting $S^1_E$ and $\partial^-_{p,q}$, which is an $\iota\oplus \tau_{p,pq}$-invariant  by 
Lemma \ref{join} (3), where part of $A$ is indicated  by dashed blue lines in Figure 5.
  And $ A \cup_{S^1_E} D^- $ is an $\iota\oplus \tau_{p,pq}$-invariant disk capping $\partial^-_{p,q}$. 
Now $$ (N*\partial _o)\cup _{\partial_o } S^-_{p,q}\cup _{\partial^-_{p,q}} ( A \cup_{S^1_E} D^-)\subset S^2*S^3$$
is $\iota\oplus \tau_{p,pq}$ invariant,  
therefore is an $f$-equivariant embedding of $\hat S^-_{p,q}$ into $S^6$, which is smooth except on the $\partial _o, \partial^-_{p,q}$ and $S^1_E$.
Note that the action is free on each 1-manifold above, so we can modify the embedding to get an equivariant smooth embedding by  Lemma \ref{modify}.

(ii) Suppose each component of $\partial^-_{p,q}$ is unknotted in $S^3$:
Then $d=(p-1; q)$ is $p-1$ or $q$, and we may assume $d=p-1$, and $\partial^{-}_{p,q}$ is a union of $d$ torus knots of type $(1, q/d)$ by Lemma \ref{lem:toptype} (3), 
and $\tau_{p,pq}$ induces a $-2\pi/(pq)$-rotation on   $\partial_o $, and $\tau^d_{p,pq}$ induces 
a $2ds\pi/(pq)$-rotation on each component of $\partial^-_{p,q}$ for some $s$ 
by Lemma \ref{lem:trans} (3).

Let $K_0$ be one of those $d$ torus knots of type $(1, q/d)$, and $K_i=\tau_{p,pq}^i(K_0)$,  $i=0,1, ..., d-1$.
Let $\iota$ be the $\frac{2\pi}d$-rotation of $S^2$ around $N$. Let $x_0$ be a point in the equator $S^1_E$, $x_i=\iota^i(x_0)$, $i=0,1, ..., d-1$.


Note the orthogonal projection from $K_0$ to $\partial_o$ provides a $\tau^d_{p,pq}$-invariant annulus $A_0$  connecting $K_0$ to $\partial _o$, and in coordinate $A_0$ is  given by $$f(r, \theta)=( \cos\frac {\pi r} 4 e^{i\theta},  \, \sin \frac{ \pi r}4 e^{iq\theta/d }),\,\, \theta\in [0, 2\pi], \, r\in [0,1]$$ 
where $f(0, \theta)=(e^{i\theta}, 0)$ gives $\partial_o$ and   $f(1, \theta)=(\frac {\sqrt 2}2 e^{i\theta}, \frac {\sqrt 2}2 e^{iq\theta/d})$ gives $K_0$. Let 
$$\partial_{o,0} =\bigcup_{y\in\partial_o }[x_0, y]_{1/2}.$$

Let  $D_0$ be the sub-disk of the geodesic disk $x_0*\partial_o$ bounded by $\partial _{o,0}.$ 
Lift  $A_0$ to $ A_0^*\subset x_0*S^3$ via the diffeomorphism
$$f(r,\theta) \mapsto [x_0, f(r,\theta)]_{1/2+r/2},$$ 
then $ A^*_0$ is a smooth annulus connecting $K_0$ and $\partial_{o,0} $. Where $K_0$, $\partial_{o}$, $A_0$, $\partial_{o,0}$, $D_0$ and $A_0^*$ are indicated in  Figure \ref{fig:CurveSur2}. 
Then $$D^*_0=D_0 \cup_{\partial _{o,0} } A^*_0$$
is an $\text{id}\oplus \tau^d_{p,pq}$-invariant disk capping on $K_0$ in $x_0*S^3\subset S^1*S^3\subset S^2*S^3$.
In Figure 6, the elements $x_0, \partial_{o,0} , A_0^*, K_0$, $D_0$ used to construct $D^*_0$ are  colored in blue.

\begin{figure}[h]
\includegraphics[width=360pt, height=180pt]{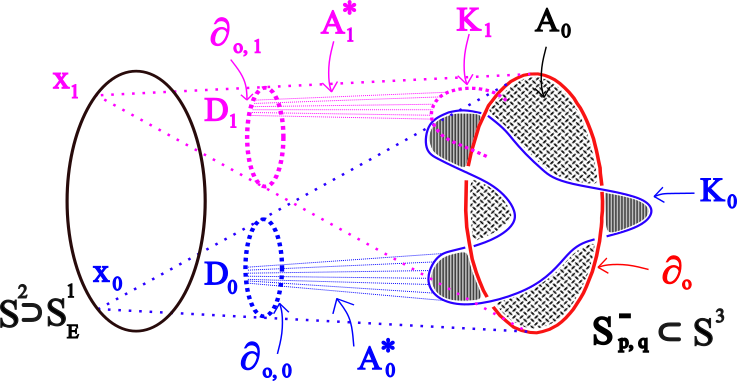}
\caption{Sketch for $f^-_{p,q}$-equivariant embedding $\hat{S}^-_{p,q}\to S^2*S^3=S^6$ when each component of $\partial^-_{p,q}$ is unknotted}\label{fig:CurveSur2}
\end{figure}

Let $\Psi =\iota\oplus \tau_{p,pq}$.
Now $D^*_i=\Psi^i(D^*_0)$ is an $id\oplus \tau^d_{p,pq}$-invariant disk capping on $K_i$ in $x_i* S^3\subset S^1*S^3\subset S^2*S^3$, $i=0,1,..., d-1$.
In Figure 6, the elements $x_1, \partial_{o,1}, A_1^*, K_1, D_1$, which are the $\Psi$ image of  $x_0, \partial_{o,0} , A_0^*, K_0, D_0$, are  colored in pink.  
Then $$(N*\partial _o)\cup _{\partial_o } S^-_{p,q}\cup_{\partial^-_{p,q}} (\cup_{i=0}^{d-1} D_i)\subset S^2*S^3 $$ is $\Psi$ invariant,  
therefore
is an $f$-equivariant embedding of $\hat S^-_{p,q}$ into $S^6$, which is smooth except on the $\partial _o, \partial^-_{p,q}$ and $\cup _{i=0}^{d-1}\Psi^i(\partial_
{o,0})$.
Note that the action is free on each 1-manifold above, so we can modify the embedding to get an equivariant smooth embedding by  Lemma \ref{modify}.

(b) The differences between the proof of $S_{p,q}$ and that of $S^\pm_{p,q}$ are: (1) there is no $\partial_o $ component; (2) $\tau_{p,q}$ acts trivially on the set of components of $\partial_{p,q}$, and it induces a $2\pi/[p\,;q]$-rotation on each component of $\partial_{p,q}$. In the case when each component of $\partial_{p,q}$ is unknotted, it is easy to see that from the proof of (ii) in (a),
$$S_{p,q}\cup_{\partial_{p,q}} (\cup_{i=0}^{d-1} D_i)\subset S^1*S^3 $$ is ${\rm id}\oplus \tau_{p,q}$ invariant, where ${\rm id}$ is the identity map on $S^1$, therefore
is an $f$-equivariant embedding of $\hat S_{p,q}$ into $S^5=S^1*S^3$. Furthermore  if $(p,q)=2$, then $\partial_{p,q}$ has two unknotted components, and we can use $S^0$ to cap those 
two components, therefore get an $f$-equivariant embedding of $\hat S_{p,q}$ into  $S^4=S^0*S^3$.
\end{proof}

 \begin{lemma}\label{modify} Let $\mathcal C$ be  a  union of finitely many disjoint simple closed curves on $\Sigma_g$.
 Suppose $\tau$ is an orientation-preserving smooth $Z_n$ action on $(\Sigma_g, \mathcal C)$ and the action
 is free on $\mathcal C$.
 If there is a $\tau$-equivariant  embedding $e: \Sigma_g\to S^m$, which is smooth except on $\mathcal C$,
 then there is a smooth $\tau$-equivariant  embedding $ \Sigma_g\to S^m$.
 \end{lemma}
 
 \begin{proof}
 We can modify it to obtain a smooth one as below: Let $\tilde \tau: S^m\to S^m$ extends $\tau$.  
 Let $\mathcal N$ be a small $\tilde \tau$-invariant regular neighborhood of $\mathcal C$ so that the action of $\tilde \tau$ on $\mathcal N$ is free,  and $\mathcal A=\mathcal N\cap e(\Sigma_g)$ 
 are $\tilde \tau$-invariant annuli with centerlines $\mathcal C$.  
 Consider the four quotients 
 $$\{N, O, A, C\}=\{\mathcal N/\tilde \tau, e(\Sigma_g)/\tilde \tau, \mathcal A/\tilde \tau, \mathcal C/\tilde \tau\}.$$
 
 Since the action of 
 ${\tilde \tau}$ on $\mathcal N$ (as well as on $\mathcal A$ and $\mathcal C$) is free, and $\mathcal N$ is a smooth manifold, $N$ is a smooth manifold;
 moreover, $A$ is a union of proper annuli in $N$ which is smooth except on their centerlines $C$. It is a classical fact that 
 we can  modify $A$ to obtain a smooth $ A'$ relative to $\partial A$, which is supported in a small neighborhood $ N'$ of $C$ in the interior of $N$. The preimage  of $(O \setminus A) \cup  A'$ gives a ${\tilde \tau}$-equivariant smooth embedding of $\Sigma_g$ into $S^m$. 
 \end{proof}
 
 \begin{proposition}\label{upper-bound1} 
 Let  $f$ be a periodic map  on $\Sigma_g$ with order $n\geq 3g$, $g>0$. 
Then $D_g(f)\le 6$. 
Furthermore, $D_1(f)\le 5$ if $f$  has order $3$;  $D_1(f)\le 4$  if $f$  has order $4$; and $D_2(f)\le 4$  if $f$  has order $6$.
\end{proposition}

\begin{table}[h]
\caption{\label{tab:3g}}
\centerline{
\begin{tabular}{|c|c|c|c|clclclcl} 
  \hline Genus & Order  & Orbifold & Model  in $S^3$ & $D_g(f)\le$ \\
  \hline $h$    & $4h+2$ & $(0:4h+2,2h+1,2)$ & $f_{2,2h+1}$, $S_{2,2h+1}$ & $6$ \\
  \hline $h$    & $4h$   & $(0:4h,4h,2)$     & $f^-_{2h,2}$,  $S^-_{2h,2}$ & $6$\\
  \hline $3h$   & $9h+3$ & $(0:9h+3,3h+1,3)$ & $f_{3,3h+1}$,  $S_{3,3h+1}$ & $6$ \\
  \hline $3h+1$ & $9h+6$ & $(0:9h+6,3h+2,3)$ & $f_{3,3h+2}$,  $S_{3,3h+2}$ & $6$\\
  \hline $3h$   & $9h$   & $(0:9h,9h,3)$     & $f^+_{3h,3}$,  $S^+_{3h,3}$ & $6$\\
  \hline $3h+1$ & $9h+3$ & $(0:9h+3,9h+3,3)$ & $f^+_{3h+1,3}$,  $S^+_{3h+1,3}$ & $6$ \\
  \hline $3h+2$ & $9h+6$ & $(0:9h+6,9h+6,3)$ & $f^-_{3h+2,3}$,  $S^-_{3h+2,3}$ & $6$ \\
  \hline $6$    & $20$   & $(0:20,5,4)$      & $f_{4,5}$,  $S_{4,5}$,& $6$\\
  \hline $9$    & $28$   & $(0:28,7,4)$      & $f_{4,7}$, $S_{4,7}$ & $6$ \\
  \hline $10$   & $30$   & $(0:30,6,5)$      & $f_{5,6}$,  $S_{5,6}$ & $6$ \\
  \hline $12$   & $36$   & $(0:36,9,4)$      & $f_{4,9}$,  $S_{4,9}$ & $6$ \\
  \hline $4$    & $12$   & $(0:12,6,4)$      & $f^-_{3,4}$,   $S^-_{3,4}$ & $6$ \\
  \hline $1$ & $3$  & $(0:3,3,3)$   & $f_{3,3}$,  $S_{3,3}$ & $5$ \\
  \hline $1$ & $4$  & $(0:4,4,2)$   & $f_{2,4}$,   $S_{2,4}$ & $4$ \\
   \hline $2$ & $6$  & $(0:6,6,3)$   & $f_{2,6}$,  $S_{2,6}$ & $4$ \\
\hline
\end{tabular}}
\end{table}

\begin{proof} Up to conjugations and powers, a periodic map $f$ on $\Sigma_g$ are determined by its order and the type of $\Sigma_g/f$ if 
$0<g\leq 3$ \cite{Br}, or 
$f$ has order $n\geq 3g$, $g>2$ \cite{Hi}. 

It is known that $D_g(f)=3$ for each free periodic map $f$ on the torus $\Sigma_1$ \cite{Ni}.
 It is also known that $D_2(f)=3$ for the periodic map $f$ of order $6$ on $\Sigma_2$ with orbifold $(0:2,2,3,3)$ by \cite[Theorem 1.5]{NWW}.

For the all remaining maps with $n\geq 3g$, we list their information  in Table 1, 
where the first three columns are genura, orders and orbifolds, which can be found in \cite[Table 1]{Br}, \cite[Theorem 3.1 and Table 1]{Hi}.

Let  $(f, \hat S)$ be either
$(f_{p,q}, \hat{S}_{p,q})$, or $(f^+_{p,q}, \hat{S}^+_{p,q})$, or $(f^-_{p,q}, \hat{S}^-_{p,q})$. 
We find a model $(f, \hat S)$ to realize the  genus $g$,  the order $n$, and  the orbifold for each row in Table 1 by Lemma \ref {lem:types}, and those $f$ and $S$ are listed in the fourth column. By Lemma  \ref{lem:toptype} one can check that  in the fourth column,
 either $\partial_{p,q}$ (resp. $\partial^{\pm}_{p,q}$) is connected or each component of $\partial_{p,q}$   (resp. $\partial^{\pm}_{p,q}$) is unknotted in $S^3$. Then by Proposition \ref{proof} (a), 
 there is an $f$-equivariant embedding  $e:\hat S\rightarrow S^6$. Hence $D_g(f)\leq 6$.
In the last three rows, each component of $\partial_{p,q}$ is unknotted, so we have $D_g(f)\leq 5$ by Proposition \ref{proof} (b), and 
in the last two  rows $(p,q)=2$, so we have $D_g(f)\leq 4$  by Proposition \ref{proof} (b) again.
All those upper bounds for $D_g(f)$ are also listed in the last column of Table 1.
\end{proof}

\begin{remark}\label{order 3} The model in Table 1 is not unique in general. For example, besides $(f_{3,3}, S_{3,3})$,  $(f^+_{1,3}, S^+_{1,3})$ is 
another model for the order 3 map on torus.
\end{remark}

\subsection{A construction from the viewpoint of orbifold}\label{subsec:orbifold}

In this section we will construct embeddings for Theorem \ref{main2}. The constructions for even $m$ are more sophisticated than 
those for odd $m$.

For a positive integer $n$, the map $\tau_{1,n}$ on $S^3$, defined in Section~\ref{subsec:Borsur},  determines a cyclic branched covering
$\pi: S^3\to \mathcal{O}_n$ of degree $n$, where $\mathcal{O}_n$ is the orientable closed 3-orbifold $S^3/\tau_{1,n}$ whose  underlying space $|\mathcal{O}_n|$ is also a 3-sphere. Let $K_n$ be the singular set of $\mathcal{O}_n$. Then $K_n$ is a trivial knot in $|\mathcal{O}_n|$ with index $n$.

\begin{figure}[h]
\includegraphics{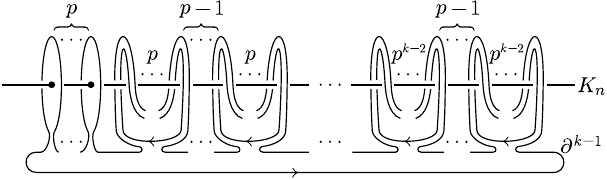}
\caption{The surface ${F}_{p^k}$ in $|\mathcal{O}_n|$ which is a band connected sum}\label{fig:BandSur}
\end{figure}

Now we let $n=p^k$, where $p$ is a prime number and $k>1$ is an integer. Let ${F}_{p^k}$ denote the compact surface in $|\mathcal{O}_{p^k}|$ as shown in Figure~\ref{fig:BandSur}. It can be obtained as a band connected sum of $p$ disks and $(k-2)(p-1)$ annuli $A_{r,s}$ where $1\leq r \leq k-2$ and $1\leq s \leq p-1$. Each disk intersects $K_n$ orthogonally at $1$ point, and $A_{r,s}$ wraps $p^r$ times around $K_n$. Let $\partial^r_s$ be the component of $\partial {F}_{p^k}$ belonging to $A_{r,s}$ and
$\partial^r=\cup_{s=1}^{p-1} \partial _s^r$, $r=1, ... , k-2$. Let $\partial^{k-1}$ be the other component, then it wraps $p^{k-1}$ times around $K_n$.

Let $ \tilde F_{p^k}=\pi^{-1}( F_{p^k})$,  the preimage of $ {F}_{p^k}$ in $S^3$. Since ${F}_{p^k}$ intersects $K_n$ transversely, $\tilde F_{p^k}$ is a connected bordered surface. Let $\hat{ F}_{p^k}$ denote the closed surface obtained by capping $\partial \tilde F_{p^k}$ with disks, then the restriction of $\tau_{1,p^k}$ on $\tilde  F_{p^k}$ extends to a periodic map on $\hat{ F}_{p^k}$. Denote this extension by $f_{p^k}$.

Denote $\pi^{-1}(\partial^r_s)$ by $\tilde \partial^r_s$, the preimage of $\partial^r_s$ in $S^3$. Since $\partial^r_s$ is a $(p^r,1)$-torus knot in $|\mathcal{O}_{p^k}|$, $\tilde \partial^r_s$ is a $(p^r,p^k)$-torus link, which is a union of $p^r$ torus knots of type $(1, p^{k-r})$. The map 
$\tau_{1,p^k}$ acts transitively on the set of the $p^r$ components of $\tilde \partial^r_s$. As in Lemma~\ref{lem:trans}, the torus link can be parameterized so that on each of its component, $\tau_{1,p^{k-r}}=(\tau_{1,p^{k}})^{p^r}$ induces a $2\pi/p^{k-r}$-rotation. Similarly, denote $\pi^{-1}(\partial ^{k-1})$ by $\tilde \partial ^{k-1}$, the preimage of $\partial^{k-1}$ in $S^3$. It is a $(p^{k-1},p^k)$-torus link, which is a union of 
$p^{k-1}$ torus knots of type $(1, p)$, where $\tau_{1,p^k}$ acts transitively on the set of components and on each component $\tau_{1,p}$ induces a $-2\pi/p$-rotation. Note that those parameterizations of torus links should induce the orientations on $\partial^r_s$ and $\partial^{k-1}$ as shown in Figure~\ref{fig:BandSur}.

Now  $\partial \tilde F_{p^k}$ consists of $k-1$ sets  $\tilde \partial ^{r}$, $r=1,..., k-1$,
where $\tilde \partial ^{k-1}$ is given and for  $r=1,..., k-2$,
$\tilde \partial ^{r}=\cup_{s=1}^{p-1} \tilde \partial^r_s.$ So $\tilde \partial ^{r}$ has $(p-1)p^r$ components.

For a fixed $r$, as described above, $\tau_{1,n}$ acts in the same manner on each subset $\tilde \partial^r_s$ of $\tilde \partial ^{r}$, $s=1, ... , p-1$
and $\tilde \partial^r_s$ consists of $p^r$ trivial knots of type $(1, p^{k-r})$. 

Now  for each $r=1, 2, ... , k-2$, the action $\tau_{1,p^{k}}$ on the set of components of $\tilde \partial^r$ has $p-1$ orbits of length $p^{r}$. We can cap $\tilde \partial ^{r}\subset S^3$ by disks in $S^3*S^1$ similar as in the proof of Proposition \ref{proof} (b), with a refinement as below. Consider the points $$\Delta^r=\{e^{\frac{2l\pi i}{p^r(p-1)}}, l=1, 2, ... ,p^r(p-1)\}\subset S^1,$$
 which are $(p-1)p^r$ points evenly distributed on the unit circle.
 Recall that $\tau_n$ is the $2\pi/n$-rotation on $S^1$. 
The action of $\tau_{p^r}$ decomposes  $\Delta^r$ into $p-1$ orbits $\Delta_s^{r}$ of length $p^r$, $s=1,..., p-1$.
Now cap the $p^r$ knots of type $(1, p^{k-r})$ in $\tilde \partial ^r_s$ from the $p^r$ points in $\Delta_s^r$, $s=1,..., p-1$,  such that 
$$\partial  \tilde F_{p^k}\cup \{\text{those $p^r(p-1)$ disks capping $\tilde \partial ^r$}\}\subset S^3*S^1$$  is invariant under the action 
$ \tau_{1,n}\oplus \tau_{p^r}$. Similarly cap the $p^{k-1}$ knots of type $(1, p)$ in $\tilde \partial ^{k-1}$ from $p^{k-1}$ evenly distribute points in $S^1$.
So by joining $S^3$ with $S^1$ $k-1$ times, we get  an  embedding 
$$e:\hat{ F}_{p^k}\rightarrow S^3\underbrace{*S^1*...*S^1}_{\text{$k-1$ times}}=S^{2k+1},$$
which is invariant under the action  
$$ \tau_{1,p^k}\oplus \tau_{p}\oplus \tau_{p^2} \oplus .... \oplus \tau_{p^{k-2}} \oplus \tau_{p^{k-1}}, \qquad (2.1)$$
therefore is an $f_{p^k}$-equivariant embedding of $\hat{ F}_{p^k}$ into $S^{2k+1}$. So we have  

\begin{lemma}\label{odd} There exists an $f_{p^k}$-equivariant embedding $e:\hat{F}_{p^k}\rightarrow S^{2k+1}$ for $k>1$.
\end{lemma}

In the proof of Lemma \ref{odd}, when $p=2$, we can replace the first $S^1$ by $S^0$ to get an equivariant embedded closed surface in $S^{2k}$. However  in this case $\tau_2$ switches the two points of $S^0$ and the map finally we get on $S^{2k}$ is orientation-reversing. So additional efforts are needed for embeddings to even dimensional spheres.

Now consider  the map $\tau_{p,p^{k+1}}$ of order $p^{k+1}$ on $S^3$ defined in Section~\ref{subsec:Borsur}, 
we have a cyclic branched covering
$\pi^*: S^3\to \mathcal{L}_{p^k}$ of degree $p^{k+1}$, where $\mathcal{L}_{p^k}=S^3/\tau_{p,p^{k+1}}$ is the orientable closed 3-orbifold, whose  underlying space is a lens space $L(p,1)$, and the singular set is a circle with index $p^k$. Denote this circle by $C_{p^k}$. 
Note that $\tau_{p,p^{k+1}}^p=\tau_{1 ,p^{k}}$ and   we have a composition
$$\pi^*=\pi'\circ \pi: S^3\to S^3/\tau_{1 ,p^{k}} =\mathcal{O}_{p^k} \to \mathcal{O}_{p^k}/\bar \tau_{p ,p^{k+1}} =\mathcal{L}_{p^k},$$
where the  map $\pi: S^3\to \mathcal{O}_{p^k}$ is given above, $\bar \tau_{p ,p^{k+1}}$ is the free action of order $p$ on $\mathcal{O}_{p^k}$ induced by $\tau_{p ,p^{k+1}}$,  and the map
$\pi': \mathcal{O}_{p^k} \to \mathcal{L}_{p^k}$ is a covering of degree $p$, which sends 
$K_{p^k}$ to $C_{p^k}$.

\begin{figure}[h]
\includegraphics{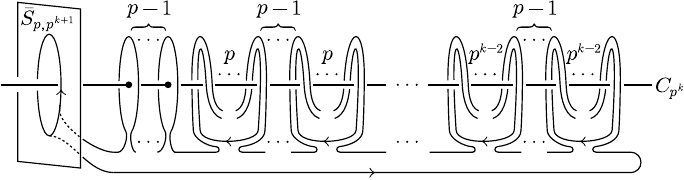}
\caption{Connect  $\bar{S}_{p,p^{k+1}}\setminus D$ to a compact surface $F$ by a band}\label{fig:BandSurPlus}
\end{figure}

We have a corresponding surface $S_{p,p^{k+1}}\subset S^3$, defined in Section~\ref{subsec:Borsur},  which is invariant under the action $\tau_{p,p^{k+1}}$. Let $\bar{S}_{p,p^{k+1}}$ denote the underlying space of $S_{p,p^{k+1}}/\tau_{p,p^{k+1}}\subset  \mathcal{L}_{p^{k}}$ with the induced orientation.
Since $(p; p^{k+1})=p$,  by Corollary ~\ref{later}, $\bar{S}_{p,p^{k+1}}$ is a $S^2$ with $p$ open disks removed and has only one singular point
(of index of $p^k$), which is the intersection of $\bar{S}_{p,p^{k+1}}$ and $C_{p^k}$. 

By removing a small open disk $D$ in $\bar{S}_{p,p^{k+1}}$ centered at the intersection,  we have the compact 
oriented surface $\bar{S}_{p,p^{k+1}}\setminus D$, then we can do the band connected sum of $\bar{S}_{p,p^{k+1}}\setminus D$ from the new boundary $\partial D$ to the oriented compact surface $F$ to get an oriented compact surface, denoted by $F^*_{p^{k+1}}$,  as shown in Figure~\ref{fig:BandSurPlus}.
Locally,   $F$ can be viewed as  the $F_{p^{k}}$ in Figure \ref{fig:BandSur} with the most left disk removed.
Now we have 
$$\partial F^*_{p^{k+1}}=\partial^1\cup \partial^2 \cup ... \cup \partial^{k-2} \cup \partial^{k-1} \cup \partial^{k+1},$$
where $\partial^r=\cup_{i=1}^{p-1} \partial _s^r$, $r=1, ... , k-2$, and $\partial ^{k-1}$ are given as before, and $\partial ^{k+1}=\partial \bar{S}_{p,p^{k+1}}$. Note that $\cup_{r=1}^{k-1} \partial ^r\subset B^3$ for some 3-ball $B^3$ around the intersection of $\bar{S}_{p,p^{k+1}}$ and $C_{p^k}$
Let $\tilde F^*_{p^{k+1}}={\pi^*}^{-1}( F^*_{p^{k+1}})$,  the preimage of ${F}^*_{p^{k+1}}$ in $S^3$, which is a connected bordered surface. Let $\hat{ F}^*_{p^{k+1}}$ denote the closed surface obtained by capping $\partial \tilde F^*_{p^{k+1}}$ with disks, then the restriction of $\tau_{p,p^{k+1}}$ on $ \tilde F^*_{p^{k+1}}$ extends to a periodic map on $\hat{ F}^*_{p^{k+1}}$. Denote this extension by $f^*_{p^{k+1}}$.

Denote $\pi^{*-1}(\partial^r_s)$ by $ \tilde \partial^r_s$, the preimage of $\partial^r_s$ in $S^3$. Then 
$$ \tilde \partial^r_s=\pi^{*-1}(\partial^r_s)=(\pi'\circ \pi)^{-1}(\partial^r_s)=\pi^{-1}\ (\pi'^{-1}(\partial^r_s)).$$
Since $\pi': \mathcal{O}_{p^k} \to \mathcal{L}_{p^k}$ is a covering of degree $p$ which sends
$K_{p^k}$ to $C_{p^k}$, and $\cup_{r=1}^{k-1} \partial ^r\subset B^3$ for some 3-ball around $C_{p^k}$,
$\pi'^{-1}((B^3, \cup_{r=1}^{k-1} \partial ^r))$ is a union of $p$ disjoint copies of $(B^3, \cup_{r=1}^{k-1} \partial ^r)_j$ around $K_{p^k}$, which is invariant
under the deck transformation.
So   the preimage of each $\partial^r_s$ in $\mathcal{O}_{p^k}$ is a union of $p$ torus knots $\partial^r_{s,j}$ of type $(1, p^r)$ in those $p$ disjoint 3-balls.
As we see before, for $\pi: S^3\to \mathcal{O}_{p^k}$,   the  preimage of each $\partial^r_{s,j}$  in $S^3$ is a torus link of type $(p^k,p^r)$.
So  $\tilde \partial^r_s$ is a union of $p$ torus links of type $(p^k, p^r)$, which gives $p^{r+1}$ torus knots of type $( p^{k-r}, 1)$. The map $\tau_{p,p^{k+1}}$ acts transitively on the set of the $p^{r+1}$ components of $\tilde \partial^r_s$ and on each of its component, 
$\tau_{1,p^{k-r}}=(\tau_{p ,p^{k+1}})^{p^{r+1}}$ induces a $2\pi/p^{k-r}$-rotation, where $1\le r \le k-2$. 
Let $\tilde \partial^{r}=\pi^{*-1}(\partial ^{r})$.  Then $\tilde \partial ^{r}=\cup_{s=1}^{p-1} \tilde \partial^r_s$, which has 
$(p-1)p^{r+1}$ components. Also $\tau_{p,p^{k+1}}$ acts transitively on the set of the $p^{k}$ components of $\tilde \partial^{k-1}$,  and  induces a $2\pi/p$-rotation on each of its component.

Finally let $\tilde \partial^{k+1}=\pi^{*-1}(\partial^{k+1})$, then $\tilde \partial^{k+1}=\partial S_{p, p^{k+1}}$ is a union of $p$ 
torus knots of type $(1, p^k)$ and  
$\tau_{p,p^{k+1}}$ keeps each component invariant. The map induces an order $p^{k+1}$ rotation on each component.
$$\partial \tilde  F^*_{p^{k+1}}=\tilde \partial^1\cup \tilde \partial^2 \cup ... \cup \tilde \partial^{k-2} \cup \tilde \partial^{k-1} \cup \tilde \partial^{k+1}.$$

Note that once more,  the action $\tau_{p,p^{k+1}}$ on the set  of components of $\tilde \partial^r$ has $(p-1)$ orbits of length $p^{r+1}$
 for $r=1, 2, ... , k-2,$, one orbit of length $p^k$ for $p=k-1$, and $p$ orbits of length 1 for $r=k+1$. 
So by joining $S^3$ with $S^1$ $k$ times as before, we get  an  embedding 
$$e:\hat{ F}^*_{p^{k+1}}\rightarrow S^3\underbrace{*S^1*...*S^1}_{\text{$k$ times}}=S^{2k+3},$$
which is invariant under the action  
$$ \tau_{p,p^{k+1}}\oplus \tau_{p^2} \oplus .... \oplus \tau_{p^{k-2}} \oplus \tau_{p^{k-1}} \oplus \tau_{p^k}\oplus {\text{id}}, \qquad (2.2)$$
therefore is an $f^*_{p^{k+1}}$-equivariant embedding of $\hat{ F}^*_{p^{k+1}}$ into $S^{2k+3}$.

Note that in (2.1)  the group action on each $S^1$ is a rotation of degree $>1$, 
and in  (2.2) the group action on the last $S^1$ is the identity. With this key difference,
when $p=2$ we can replace the last $S^1$ by $S^0$ to get an equivariant embedded closed surface in $S^{2k+2}$, still in the orientable category. 

Finally note that $\hat{F}^*_{p^{k+1}}/f^*_{p^{k+1}}$ is a sphere having $p$ points of index $p^{k+1}$, $p-1$ points of index $p^{k+1-r}$ for $1\leq r<k$, and one point of index $p$: Refer to Figure \ref{fig:BandSurPlus}, it is easy to see that the underlying space of $\hat{F}^*_{p^{k+1}}/f^*_{p^{k+1}}$ is a sphere. Then, $p$ points of index $p^{k+1}$ correspond to the $p$ boundary components of $\bar S_{p, p^{k+1}}$;
$p-1$ points of index $p^{k}$ correspond to the centers of those $p-1$ disks; $p-1$  point of index $p^{k-r}$  correspond those $p-1$ $\partial^r_s$ for $1\leq r<k-1$; a  point of index $p$ corresponds to $\partial^{k-1}$. 

\begin{proposition}\label{even}  (1) For $k>2$,
there exists an $f^*_{p^{k}}$-equivariant embedding $e:\hat F^*_{p^{k}}\rightarrow S^{2k+1}$. Moreover, when $p=2$, there is an $f^*_{2^{k}}$-equivariant embedding $e:\hat F^*_{2^{k}}\rightarrow S^{2k}$.

 (2)  $\hat{F}^*_{p^k}/f^*_{p^k}$ is a sphere having $p$ points of index $p^k$, $p-1$ points of index $p^{k-r}$ for $1\leq r<k-1$, and one point of index $p$. 
 
 (3) The genus of $\hat{F}^*_{p^k}$ is $\frac{(k-1)(p-1)p^k}{2}-p^{k-1}+1.$
\end{proposition}

\begin{proof} (1) and (2) follow the last three paragraphs before Proposition \ref{even},  by changing $k+1$ to $k$ and $k>1$ to $k>2$.

(3) follows (2) and that the degree of the cyclic branch covering $\hat{F}^*_{p^k}\to \hat{F}^*_{p^k}/f^*_{p^k}$ is $p^k$:
Suppose  the orientable closed surface $\hat{\mathcal F}_{p^{k+1}}$ has genus $g$.
By the Riemann--Hurwitz formula,
\[2-2g=p^k(2-\sum_{r=1}^{k-2}(p-1)(1-\frac{1}{p^{k-r}})-(1-\frac{1}{p})-p(1-\frac{1}{p^k})).\]
Hence \[g=\frac{(k-1)(p-1)p^k}{2}-p^{k-1}+1.\]

\end{proof}




\section{Proofs of the theorems}\label{sec:proof}
Suppose that there is an $f$-equivariant embedding $e:\Sigma_g\rightarrow S^m$. By definition we have an orientation-preserving periodic map $\tilde{f}$ on $S^m$ such that $\tilde{f}$ has the same order as $f$ and $e\circ f=\tilde{f}\circ e$. In this section we will analyze the fixed point sets of the powers of $\tilde{f}$ and their tangent spaces by using some classical facts in differential and algebraic topology and differential geometry, in particular the  Smith theory, to give the required lower bounds of $D_g(f)$. Let $Fix(\cdot)$ denote the fixed point set of a periodic map. We will use the following properties about the fixed point set.

\begin{proposition}\label{lem:fix}
Let $\tau$ be a smooth periodic map (isometry) on a  smooth (closed Riemannian) manifold $M$. Then

(1) $Fix(\tau)$ is a smooth (totally geodesic closed) submanifold of $M$.

(2) When $M$ is orientable and $\tau$ is orientation-preserving, the codimension of each component of $Fix(\tau)$ is even. 


(3) If $M$ is $S^m$ and the order of $\tau$ is a power of a prime number $p$, then $Fix(\tau)$ is a $\mathbb{Z}_p$-homology sphere, where $\mathbb{Z}_p=\mathbb{Z}/p\mathbb{Z}$.

(4) Let $M$ be a smooth rational homology sphere. If two closed oriented submanifolds of $M$ meet transversely, then their
algebraic intersection number is 0.
\end{proposition}

For (1), see \cite{Wal} and \cite[Theorem 5.1]{Kob}, (2) can be proved by choosing an equivariant Riemannian metric on $M$ and considering the induced map of $\tau$ on the tangent space of a fixed point, and (3) follows from  Smith theory \cite{Sm}.
(4) is a standard fact in differential topology. Note that in (1) $Fix(\tau)$ need not be connected, and in (3) $Fix(\tau)$ is disconnected only if it is a set of two points. This is a key ingredient in our proofs.

By Proposition~\ref{lem:fix}, if $f$ has order $n=p^rq$, where $p$ is a prime number and $r\geq 1$, then $Fix(\tilde{f}^q)$ is a smooth $\mathbb{Z}_p$-homology sphere and $m-\dim Fix(\tilde{f}^q)$ is even.

\begin{lemma}\label{>3}
If $f$ has order $n>1$ and $\Sigma_g/f$ has type $(\bar{g}:n,a,b)$, $a, b>1$, then $D_g(f)\ge 4$.
\end{lemma}

\begin{proof}
Suppose that $f$ extends to $\tilde{f}$ on $S^m$ with respect to $e$. 
If $m=3$,
then $|S^3/\tilde f|$ is a rational homology 3-sphere, and   the singular set of $S^3/\tilde f$ is  a closed 1-manifold (since $\tilde f$ is orientation preserving), which meet the closed orientable surface $|\Sigma_g/f|$ transversely in three points, the singular points of $\Sigma_g/f$, so their algebraic intersection number in nonzero, which contradicts Proposition \ref{lem:fix} (4).
\end{proof}

\begin{proposition}\label{ind}
Suppose $f$ has order $n>1$ and $\Sigma_g/f$ has type $(\bar{g}:n,a,b)$, $a, b>1$. Then

(1) if $a,b>1$ and $(a\,;b)=1$, then $D_g(f)\geq 6$;

(2) if $n=a>2b>2$, then $D_g(f)\geq 6$;

(3) if $n=a=b$, then $D_g(f)\geq 5$.
\end{proposition}


\begin{proof}

Note that we always have $a\mid n$, $b\mid n$ and $a,b>1$. Write $n=ar=bs$. Recall that each singular point $x_i$ of index $n_i$  corresponds to an $f$-orbit of length $n/n_i$ and each point in this orbit  has stabilizer of order $n_i$.

In (1), since $(a\,;b)=1$, and $a\mid n$, $b\mid n$, $a,b>1$, we have $a,b<n$. So $f$ has exactly $1$ fixed point. Denote it by $x$. If a nontrivial power $f^t$ has a fixed point other than $x$, then either $r\mid t$ or $s\mid t$, and the fixed point is a preimage of the point in $\Sigma_g/f$ of index $a$ or $b$, respectively. Choose a prime factor $p$ of $a$ and a prime factor $q$ of $b$. Write $n=pu=qv$. Then, since $r\mid u$ and $s\nmid u$, $f^u$ has exactly $r+1$ fixed points, corresponding to the points of indices $n$ and $a$. Similarly $f^v$ has exactly $s+1$ fixed points, corresponding to the points of indices $n$ and $b$. Note that $r+1,s+1\geq 3$. Hence $Fix(\tilde{f}^u)$ and $Fix(\tilde{f}^v)$, as homology spheres  by Proposition \ref{lem:fix} (3), furthermore, have dimension at least $1$. They have a common point $e(x)$, which is fixed by $\tilde{f}$, and they are not contained in each other.

Let $V$ be the tangent space of $S^m$ at $e(x)$, and let $\tilde{f}_\ast$ be the map on $V$ induced by $\tilde{f}$. We view $V$ as a vector space over $\mathbb{C}$. Then, since $\tilde{f}$ has order $n$, $\tilde{f}_\ast^n$ must be the identity, and there is a decomposition of $V$ into invariant subspaces as
\[V=\bigoplus_{d\mid n}V(d),\]
where $V(d)$ is a direct sum of the eigenspaces with eigenvalues the primitive $d$-th roots of unity. Since $\tilde{f}$ is orientation-preserving, $\dim V(d)$ is even for $d>1$.

Now let $V_1$ and $V_2$ be the tangent spaces of $Fix(\tilde{f}^u)$ and $Fix(\tilde{f}^v)$ at $e(x)$, respectively. By choosing an $\tilde{f}$-invariant Riemannian metric on $S^m$, the exponential map at $e(x)$ will provide a correspondence between $Fix(\tilde{f}_\ast^t)$ and the component of $Fix(\tilde{f}^t)$ containing $e(x)$. According to the above decomposition of $V$, we have
\[V_1=\bigoplus_{d_1\mid u}V(d_1), V_2=\bigoplus_{d_2\mid v}V(d_2).\]

Since $Fix(\tilde{f}^u)$ and $Fix(\tilde{f}^v)$ are connected totally geodesic submanifolds (Proposition \ref{lem:fix} (1)) and are not contained in each other, $V_1$ and $V_2$ are not contained in each other. So there exist nonzero subspaces $V(d_1)$ and $V(d_2)$ so that $d_1\mid u$, $d_1\nmid v$ while $d_2\mid v$, $d_2\nmid u$. So $d_1$, $d_2$ and $n$ are three distinct integers. On the other hand, since $\tilde f_*$ is a rotation of order $n$ on $T_{e(x)}e(\Sigma_g)$, the tangent space of $e(\Sigma_g)$ at $e(x)$,   $T_{e(x)}e(\Sigma_g)\subset V(n)$, which implies that $\dim V(n)\ge 2$. Note that $\text{dim}V(d_1)$ and $\text{dim}V(d_2)$ are even and nonzero. So
\[m=\dim V\geq \dim V(d_1)+\dim V(d_2)+\dim V(n)\geq 6.\]

In (2), since $n=a>2b$, $f$ has exactly $2$ fixed points. If a nontrivial power $f^t$ has another fixed point $y$, then $s\mid t$, and $y$ is a preimage of the point in $\Sigma_g/f$ of index $b$. Choose a prime factor $q$ of $b$ and write $n=qv$. Since $s\mid v$, $f^v$ has exactly $s+2$ fixed points.  Let $F$ be a component of $Fix(\tilde{f})$. Since $Fix({f})\subset Fix({f}^v)$ is proper, $Fix(\tilde{f})\subset Fix(\tilde{f}^v)$ is proper. Then by Proposition~\ref{lem:fix} (2)  and (3), $Fix(\tilde{f}^v)$ is a smooth homology sphere, and
\[\dim Fix(\tilde{f}^v)\geq \dim F+2\geq 2.\]


Note that the intersection of $Fix(\tilde{f}^v)$ and $e(\Sigma_g)$ contains exactly $s+2$ points. By Lemma \ref{>3}, $m\ge 4$. If $m=4$, 
then $\text{dim}Fix(\tilde f^v)=2$ by Proposition \ref{lem:fix} (2), so
$Fix(\tilde{f}^v)$ and $e(\Sigma_g)$ intersect transversely. Since $\tilde{f}$ is orientation-preserving on $S^m$ and $e(\Sigma_g)$, and it preserves $Fix(\tilde{f}^v)$, it also preserves the orientations of $Fix(\tilde{f}^v)$. So the $s$ points in the orbit of $e(y)$ have the same sign. Since $n>2b$,  $s\geq 3$. Since $s\geq 3$, the algebraic intersection number of $Fix(\tilde{f}^v)$ and $e(\Sigma_g)$ is nonzero, which is a contradiction 
by Proposition \ref{lem:fix} (4). If $m=5$, then $\dim F\geq 1$ by Proposition~\ref{lem:fix} (2), and $\dim Fix(\tilde{f}^v)= 3$. Then, by the same reason, we have a contradiction. Hence $m\geq 6$.


 In (3), since $n=a=b$, $f$ has exactly $3$ fixed points, and any nontrivial power $f^t$ has no other fixed points. Choose a prime factor $p$ of $n$ and write $n=pu$. Then $f^u$ has exactly $3$ fixed points. So $Fix(\tilde{f}^u)$ is a smooth homology sphere of dimension at least $1$. 
 By Lemma \ref{>3}, $m\ge 4$.
 If $m=4$ then $\dim Fix(\tilde{f}^u)= 2$ by Proposition~\ref{lem:fix} (2), and $Fix(\tilde{f}^u)$ and $e(\Sigma_g)$ will intersect transversely at $3$ points, so the algebraic intersection number is odd, which is impossible by Proposition~\ref{lem:fix} (4). So $m\geq 5$.
 
\end{proof}

\begin{remark}\label{ind*}
In Proposition~\ref{ind} (1), the condition ``$(a\,;b)=1$'' can be replaced by the following weaker one: there exist two integers $p$ and $q$, which are powers of prime numbers, such that $p\mid a$, $p\nmid b$ and $q\mid b$, $q\nmid a$. The proof is almost the same.
\end{remark}

\begin{lemma}\label{lem:seq}
Assume that $\Sigma_g/f$ has $k$ points of indices $p^{r_j}q_j$, $1\leq j\leq k$, where $p$ is a prime number, $p\nmid q_j$ for $1\leq j\leq k$ and $r_1>\ldots>r_k>0$. Then $D_g(f)\geq 2k$. Moreover, if $f$ has at least $3$ fixed points, then $D_g(f)\geq 2k+1$.
\end{lemma}

\begin{proof}
Let $n$ be the order of $f$. For $1\leq j\leq k$, write $n=p^{r_j}q_js_j$, let $v_j=q_js_j$, and let $F_j$ be the set of preimages of the point in $\Sigma_g/f$ of index $p^{r_j}q_j$. If a power $f^t$ fixes a point in $F_j$, then $s_j\mid t$. Since $p\nmid q_j$ for $1\leq j\leq k$ and $r_1>\ldots>r_k>0$, $s_i\mid v_j$ if and only if $i\leq j$. So $F_i\subseteq Fix(f^{v_j})$ if and only if $i\leq j$. Note that $v_i\mid v_j$ if and only if $i\leq j$. Hence we have a sequence of strictly inclusions
\[Fix(f^{v_1})\subset Fix(f^{v_2})\subset\ldots\subset Fix(f^{v_k})\subset \Sigma_g.\]
Then in $S^m$ we also have a sequence of strictly inclusions
\[Fix(\tilde{f}^{v_1})\subset Fix(\tilde{f}^{v_2})\subset\ldots\subset Fix(\tilde{f}^{v_k})\subset S^m.\]
By Lemma~\ref{lem:fix}, all $Fix(\tilde{f}^{v_j})$ are $\mathbb{Z}_p$-homology spheres with even codimension, so
\[\dim Fix(\tilde{f}^{v_1})<\dim Fix(\tilde{f}^{v_2})<\ldots<\dim Fix(\tilde{f}^{v_k})<m\]
and $m\geq\dim Fix(\tilde{f}^{v_1})+2k\geq 2k$. If $f$ has at least $3$ fixed points, then $Fix(\tilde{f}^{v_1})$ is not the set of two points. Hence $\dim Fix(\tilde{f}^{v_1})\geq 1$ and $m\geq 2k+1$.
\end{proof}


\begin{proof}[Proof of Theorem~\ref{main1}]
According to Table~\ref{tab:3g}, we have the required upper bounds of $D_g(f)$. As we mentioned in the proof of Proposition 
\ref{upper-bound1},
 the map given by the bottom line of Table~\ref{tab:3g} has $D_g(f)=3$. 
Since each $(\Sigma_g, f)$ in the remaining lines of Table~\ref{tab:3g}  has $D_g(f)>3$ by Lemma \ref{>3}, so the part of $D_g(f)=4$ is verified.
The part of $D_g(f)=5$ is verified by Proposition~\ref{ind} (3). The part $D_g(f)=6$ is verified by Proposition~\ref{ind} (1) and Proposition \ref{ind} (2), except for the map of order $12$ on $\Sigma_4$, which has orbifold type $(0:12,6,4)$. Because $3\mid 6$, $3\nmid 4$ and $2^2\mid 4$, $2^2\nmid 6$, this case can be verified by Remark~\ref{ind*}.
\end{proof}

\begin{proof}[Proof of Theorem~\ref{main2}] By Theorem \ref{main1}, we may assume that $m\ge 6$.

In Section~\ref{subsec:orbifold}, for $k>2$, we have constructed the maps $f^*_{p^k}$ 
such that $D_g(f^*_{p^k})\leq 2k+1$, and $D_g(f^*_{2^k})\leq 2k$  by Proposition \ref{even}.

By Proposition \ref{even}, $\hat{F}^*_{p^k}/f^*_{p^k}$ is a sphere having $p$ points of index $p^k$, $p-1$ points of index $p^{k-r}$ for $1\leq r<k-1$, and $1$ point of index $p$. So by Lemma~\ref{lem:seq}, $D_g(f^*_{p^k})\geq 2k$. Hence $D_g(f^*_{2^k})=2k$. If $p\geq 3$, then there are at least $3$ points of index $p^k$. So by the moreover part of Lemma \ref{lem:seq}, $D_g(f^*_{p^k})=2k+1$. 
 
 Because the condition in Lemma~\ref{lem:seq} only concerns about singular points, by applying surgeries equivariantly at preimages of a regular point, we can enlarge $g$ to obtain infinitely many required $(\Sigma_g, f)$. 
\end{proof}










Finally we restate Proposition \ref{thm:s} more concretely as Propositions \ref{thm:g} and \ref{thm:s1}.


\begin{proposition}\label{thm:g} $D_g(f)=\hat{D}_g(f)$ for all maps $f$ in the proof of Theorem \ref{main2}. 
Hence for any integer $k>1$, there exist infinitely many $(\Sigma_g, f)$, where $f$ is a periodic map on $\Sigma_g$, such that $\hat{D}_g(f)=k$.
\end{proposition}

\begin{proof}
Note that Lemma~\ref{lem:seq} still holds even if the embedding is not smooth. So Proposition~\ref{thm:g} follows from the proof of Theorem \ref{main2}.
\end{proof}

\begin{proposition}\label{thm:s1}
For  all maps $f$  with  ${D}_g(f)=6$  in Theorem \ref{main1}, except the one on line 12 of Table 1, $\hat{D}_g(f)\le 4$.
\end{proposition}

\begin{proof} 
According to  Table 1, for  all maps $f$  with  ${D}_g(f)=6$  in Theorem \ref{main1}, except the one on line 12 of Table 1,
 we have constructed an $f$-equivariant corresponding bordered surface $S\subset S^3$ so that  $\partial S$  has at most $2$ 
 components by Lemma~\ref{lem:toptype}.  Note  $S\subset S^3\subset S^4$, by making join from the north pole to one component of $\partial S$, (and from south pole to the other component of $\partial S$, if $\partial S$ has two components)
we have an $f$-equivariant embedding  from $\Sigma_g$ into $S^4$. So $\hat{D}_g(f)\le 4$. 
\end{proof}

\bibliographystyle{amsalpha}

\end{document}